\documentclass[leqno]{amsart}
\usepackage[normal]{boilerplate}

\title{Categorifying rationalization}

\author{Clark Barwick}
\address{Massachusetts Institute of Technology, Department of Mathematics, 77 Massachusetts Avenue, E17-332, Cambridge, MA 02139-4307, USA}
\email{clarkbar@gmail.com}

\author{Saul Glasman}
\address{University of Minnesota, School of Mathematics, Vincent Hall, 206 Church St. SE, Minneapolis, MN 55455, USA}
\email{saulglasman0@gmail.com}

\author{Marc Hoyois}
\address{Department of Mathematics, Massachusetts Institute of Technology, 77 Massachusetts Avenue, Cambridge, MA 02139-4307, USA}
\email{hoyois@mit.edu}

\author{Denis Nardin}
\address{Department of Mathematics, Massachusetts Institute of Technology, 77 Massachusetts Avenue, Cambridge, MA 02139-4307, USA}
\email{nardin@math.mit.edu}

\author{Jay Shah}
\address{Department of Mathematics, Massachusetts Institute of Technology, 77 Massachusetts Avenue, Cambridge, MA 02139-4307, USA}
\email{jshah@math.mit.edu}

\begin{document}

\begin{abstract} We solve a problem proposed by Khovanov by constructing, for any set of primes $S$, a triangulated category (in fact a stable $\infty$-category) whose Grothendieck group is $S^{-1}\mathbf{Z}$. More generally, for any exact $\infty$-category $E$, we construct an exact $\infty$-category $S^{-1}E$ of equivariant sheaves on the Cantor space with respect to an action of a dense subgroup of the circle. We show that this $\infty$-category is precisely the result of categorifying division by the primes in $S$. In particular, $K_n(S^{-1}E)\cong S^{-1}K_n(E)$.
\end{abstract}

\maketitle

It is a peculiar fact that rationalized algebraic $K$-groups have largely remained out of reach of algebraic techniques. For example, the rationalized $K$-groups of a number field $F$ were computed by Borel \cite{MR0387496}: for $n\geq 2$,
\[
\dim K_n(F)\otimes\QQ=\begin{cases} 0 & \text{if }n\equiv 0\mod 2;\\ r_1+r_2 & \text{if }n\equiv 1\mod 4;\\ r_2 & \text{if }n\equiv 3\mod 4,
\end{cases}
\]
where $r_1$ is the number of real places and $r_2$ is the number of complex places of $F$. But Borel's proof depends upon a delicate analysis of invariant differential forms on the Borel--Serre compactification of a symmetric space. As far as we know, no algebraic approach to this computation has appeared in the literature.

For function fields, the situation is at least as dire. For example, we have the following.
\begin{cnj*}[Parshin] If $X$ is a smooth projective variety over a finite field, then $K_n(X)\otimes\QQ=0$ for any $n\geq 1$.
\end{cnj*}
\noindent But only when the dimension of $X$ is $0$ or $1$ is this assertion known.

The task of this paper is to \emph{categorify rationalization}, in order to get a more explicit grasp on rational $K$-theory classes. That is, we introduce explicit categories of \emph{divisible objects} whose $K$-theory gives the rational $K$-theory directly.

More precisely, if $S$ is a set of prime numbers, then for any exact $\infty$-category $E$ (in particular, for any exact ordinary category or any stable $\infty$-category \cite{MR3427577}), we construct here an exact $\infty$-category $S^{-1}E$ such that $K(S^{-1}E)\simeq S^{-1}K(E)$ as spectra, and, consequently,
\[
K_{\ast}(S^{-1}E)\cong S^{-1}K_{\ast}(E)
\]
as graded abelian groups. 

When $E$ is an idempotent-complete stable $\infty$-category, we can offer an explicit -- though perhaps unwieldy -- characterization of $S^{-1}E$: it is an $\infty$-category of what we call \emph{$S$-divisible objects}. These are sequences $\{X_i\}$ of objects $X_i$ of $\Ind E$, indexed over the various products $i$ of the primes in $S$, along with suitably compatible identifications, when $m$ divides $n$, between the object $X_m$ and the $n/m$-fold direct sum $X_n\oplus X_n \oplus\cdots\oplus X_n$, all subject to a finiteness condition.

Our main theorem goes a step still further, and identifies $S^{-1}E$ as an $\infty$-category of sheaves of objects of $\Ind E$ on the Cantor space $\Omega$ that are equivariant with respect to a free action (Cnstr. \ref{prp:Q/Z-action}) of the \emph{$S$-adic circle}
\[\TT_S\coloneq S^{-1}\ZZ/\ZZ\]
on $\Omega$. When $E$ is the $\infty$-category of coherent complexes on a reasonable scheme $X$, we may think of $\Omega$ as an affine scheme with its $S$-adic circle action, and we prove:
\begin{thm*} One has an equivalence of $\infty$-categories
\[
S^{-1}\categ{IndCoh}(X)\simeq\categ{IndCoh}^{\TT_S}(X\times \Omega)
\]
between the $S$-divisible ind-coherent complexes on $X$ and the $\infty$-category of $\TT_S$-equi\-variant ind-coherent complexes on $X\times \Omega$.
\end{thm*}

\noindent We deduce that
\[
S^{-1}G_n(X)\cong G_n^{\TT_S}(X\times \Omega)\,,\text{\quad and in particular\quad}G_n(X)\otimes\QQ\cong G_n^{\QQ/\ZZ}(X\times \Omega)\,;
\]
that is, the rationalized $G$-theory of $X$ is the $\QQ/\ZZ$-equi\-variant $G$-theory of $X\times \Omega$.

This paper thus solves problems posed by Khovanov \cite[2.3 and 2.4]{Khovanov}, who sought such a ``categorification of division.'' In particular, he asked for a triangulated category whose Grothendieck group is $\QQ$, and more generally, one whose Grothendieck group is $m^{-1}\ZZ$ for an integer $m$. In fact, for any field $k$, the stable $\infty$-category $\categ{QCoh}^{m^{-1}\ZZ/\ZZ}(\Spec k\times \Omega)$ of $m^{-1}\ZZ/\ZZ$-equivariant sheaves of complexes of $k$-vector spaces on $\Omega$ is the localization of the derived category of $k$ away from $m$. The compact objects therein have not only the desired Grothendieck group $m^{-1}\ZZ$, but one even has
\[
m^{-1}K_n(k)\cong G_n^{m^{-1}\ZZ/\ZZ}(\Spec k\times \Omega).
\]
The slogan is thus: \emph{Vector spaces with rational dimension are circle-equivariant sheaves of complexes on the Cantor space.}

Finally, though our motivation was to contemplate rational algebraic $K$-theory, we must note that nowhere have we really used anything special about the functor $K$, save only that it preserves finite products and filtered colimits. Any functor with this property (e.g., topological Hochschild homology) can replace $K$ in the assertions above. This reflects the fact that our procedure really inverts the primes in $S$ at the \emph{categorical} level.

\subsection*{Acknowledgements} We thank R. Bezrukavnikov for a helpful conversation about this paper and for pointing us to Khovanov's conjecture. We thank A. Putman for an encouraging conversation.


\section{Localizations}

\begin{rec} An abelian group $E$ is \emph{$S$-local} if and only if, for product $k$ of primes in $S$, the multiplication by $k$ map $k\colon\fromto{E}{E}$ is an isomorphism.
\end{rec}

More generally, we have the following.
\begin{dfn} Suppose $C$ an $\infty$-category with direct sums. For any object $E$ of $C$, and for any natural number $k$, write $kE$ for the $k$-fold direct sum $E\oplus E\oplus \cdots\oplus E$. The composite
\[
E\to kE\to E
\]
of the codiagonal followed by the diagonal deserves the name \emph{multiplication by $k$}. We will say that $E$ is \emph{$S$-local} if and only if, for any product $k$ of primes in $S$, the multiplication by $k$ map $k\colon\fromto{E}{E}$ is an equivalence.
\end{dfn}

This recovers, e.g., the notion of $S$-locality for spectra.

\begin{ntn} Let $\Phi_S$ denote the ordinary category in which an object is a (positive) natural number that is a product of elements of $S$, and a morphism $\fromto{m}{n}$ is a natural number $k$ such that $n=mk$.
\end{ntn}

We will show in \S \ref{app:constructESandES} that every object $E$ of an $\infty$-category $C$ with direct sums determines a functor
\[E[S]\colon\fromto{\Phi_S}{C}\]
that carries every object to $E$ and every morphism $k\colon\fromto{m}{n}$ to the morphism $k\colon\fromto{E}{E}$, as well as a dual diagram
\[E[S]^{\vee}\colon\fromto{\Phi_S^{\op}}{C}\]
that carries every object to $E$ and every morphism $k\colon\fromto{m}{n}$ to the morphism $k\colon\fromto{E}{E}$.

The proof of the following is easy.
\begin{prp} Suppose $C$ an $\infty$-category that admits direct sums and filtered colimits. Then the following are equivalent for an object $E$ of $C$.
\begin{itemize}
\item The object $E$ is $S$-local.
\item The functor $E[S]$ is essentially constant.
\item The natural map $\fromto{E}{\colim E[S]}$ is an equivalence in $C$.
\end{itemize}
\end{prp}

\begin{ntn} If $C$ is an $\infty$-category that admits direct sums and filtered colimits, then we write $S^{-1}\colon\fromto{C}{C}$ for the functor $\goesto{E}{\colim E[S]}$.
\end{ntn}

\begin{wrn} It is tempting to believe that $S^{-1}\colon\fromto{C}{C}$ is a localization functor onto the full subcategory spanned by the $S$-local objects. This is true when $C$ is $\Ab$ or $\Sp$. However, it isn't true in general: see Warning~\ref{wrn:Slocal}.
In order for $S^{-1}E$ to be $S$-local, it is sufficient that for any $p\in S$, there exist $N\geq 2$ such that the cyclic permutation of $p^N\colon\fromto EE$ is homotopic to the identity.
\end{wrn}

\section{The effective Burnside $\infty$-category and the functors $E[S]$ and $E[S]^{\vee}$}\label{app:constructESandES} We give a precise construction of the functors $E[S]$ and $E[S]^{\vee}$ for any object $E$ of any $\infty$-category $C$ that admits direct sums and filtered colimits.

To this end, let $A^{\eff}(\Fin)$ denote the effective Burnside $\infty$-category of finite sets \cite{M1}. (This is in fact a $2$-category.) We have shown that this is the Lawvere theory of $E_{\infty}$ objects. That is, for any $\infty$-category $D$ with all finite products, there is an equivalence
\[\categ{CAlg}(D^{\times})\simeq\Fun^{\times}(A^{\eff}(\Fin),D),\]
where $\Fun^{\times}$ denotes the $\infty$-category of product-preserving functors. Equivalently, $A^{\eff}(\Fin)$ can be identified with the $\infty$-category of free, finitely generated $E_{\infty}$ spaces.

Now since $C$ has direct sums, every object is an $E_{\infty}$-algebra in a unique way. That is, the forgetful functor
\[
\fromto{\categ{CAlg}(C^{\times})}{C}
\]
is an equivalence. Consequently, the functor 
\[
\fromto{\Fun^{\times}(A^{\eff}(\Fin),C)}{C}
\]
given by evaluation at the one-point set $\langle1\rangle\coloneq\{0\}$ is an equivalence. Select, once and for all, a homotopy inverse $F$ to this equivalence. Now in order to construct $E[S]$ and $E[S]^{\vee}$ for any object $E$ of $S$, we need only to define a functor
\[
M_S\colon\fromto{\Phi_S}{A^{\eff}(\Fin)}
\]
that carries each natural number in $\Phi_S$ to the singleton, and every map $\fromto{m}{n}$ given by $n=mk$ to the span
\[
\langle1\rangle \ot \langle k\rangle\to \langle1\rangle,
\]
where
\[\langle k\rangle\coloneq\{0,1,\dots,k-1\}\,.\]
We then obtain $E[S]$ as the composite $F(E)\circ M_S$, and we obtain $E[S]^{\vee}$ as the composite $F(E)\circ D\circ M_S^{\op}$, where $D\colon\equivto{A^{\eff}(\Fin)^{\op}}{A^{\eff}(\Fin)}$ is the duality functor.

In fact it will be useful to define a functor
\[
\widetilde{M}_S\colon\fromto{O(\Phi_S)}{A^{\eff}(\Fin)}\,,
\]
where $O(\Phi_S)\coloneq\Fun(\Delta^1,\Phi_S)$ is the arrow category of $\Phi_S$, such that the precomposition of $\tilde{M}_S$ with the inclusion $\Phi_S\subseteq O(\Phi_S)$ sending every object to the identity on it is the required functor $M_S$.

To define $\widetilde M_S$ carefully, if $n=mk$, then we define two maps
\[
p_{m|n}\colon\fromto{\langle n\rangle}{\langle m\rangle}\text{\quad and\quad}j_{m|n}\colon\fromto{\langle n\rangle}{\langle m\rangle}
\]
by the formulas
\[
p_{m|n}(i)\coloneq\left\lfloor \frac{i}{k}\right\rfloor\text{\quad and\quad}j_{m|n}(i)\coloneq i\mod m.
\]

Now for any $p$-simplex
\[(m_0|n_0)|(m_1|n_1)|\cdots|(m_p|n_p)\]
of $O(\Phi_S)$ (by which we mean that $m_s|m_{s+1}$ and $n_t|n_{t+1}$) in which $n_t=k_{s,t}m_s$, the $p$-simplex
\[
\widetilde{M}_S((m_0|n_0)|(m_1|n_1)|\cdots|(m_p|n_p))\in A^{\eff}(\Fin)_p
\]
will be the diagram
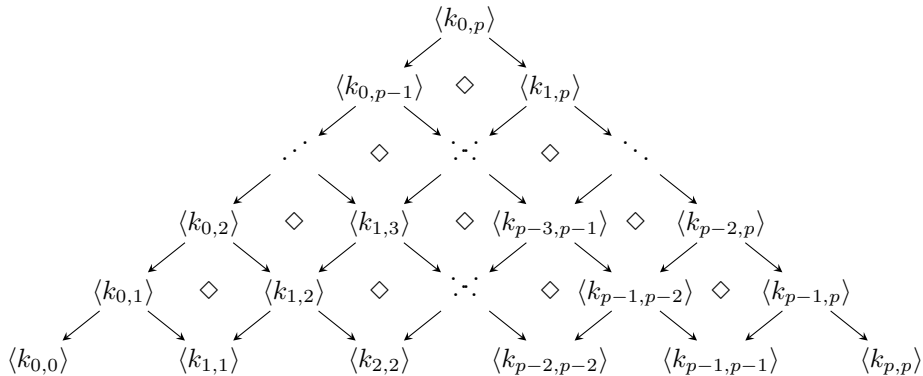
\begin{equation*}
\begin{tikzpicture} 
\matrix(m)[matrix of math nodes, 
row sep={6ex,between origins}, column sep={7.5ex,between origins}, 
text height=1.5ex, text depth=0.25ex] 
{&&&&&\langle k_{0,p}\rangle&&&&&\\
&&&&\langle k_{0,p-1}\rangle&\Diamond&\langle k_{1,p}\rangle&&&&\\
&&&\adots&\Diamond&\xdots&\Diamond&\ddots&&&\\
&&\langle k_{0,2}\rangle&\Diamond&\langle k_{1,3}\rangle&\Diamond&\langle k_{p-3,p-1}\rangle&\Diamond&\langle k_{p-2,p}\rangle&&\\
&\langle k_{0,1}\rangle&\Diamond&\langle k_{1,2}\rangle&\Diamond&\xdots&\Diamond&\langle k_{p-1,p-2}\rangle&\Diamond&\langle k_{p-1,p}\rangle&\\
\langle k_{0,0}\rangle&&\langle k_{1,1}\rangle&&\langle k_{2,2}\rangle&&\langle k_{p-2,p-2}\rangle&&\langle k_{p-1,p-1}\rangle&&\langle k_{p,p}\rangle\\}; 
\path[>=stealth,<-,font=\scriptsize] 
(m-2-5) edge (m-1-6) 
(m-2-7) edge (m-1-6)
(m-3-4) edge (m-2-5)
(m-3-6) edge (m-2-5)
(m-3-6) edge (m-2-7)
(m-3-8) edge (m-2-7)
(m-4-3) edge (m-3-4)
(m-4-5) edge (m-3-4)
(m-4-5) edge (m-3-6)
(m-4-7) edge (m-3-6)
(m-4-7) edge (m-3-8)
(m-4-9) edge (m-3-8)
(m-5-2) edge (m-4-3)
(m-5-4) edge (m-4-3)
(m-5-4) edge (m-4-5)
(m-5-6) edge (m-4-5)
(m-5-6) edge (m-4-7)
(m-5-8) edge (m-4-7)
(m-5-8) edge (m-4-9)
(m-5-10) edge (m-4-9)
(m-6-1) edge (m-5-2)
(m-6-3) edge (m-5-2)
(m-6-3) edge (m-5-4)
(m-6-5) edge (m-5-4)
(m-6-5) edge (m-5-6)
(m-6-7) edge (m-5-6)
(m-6-7) edge (m-5-8)
(m-6-9) edge (m-5-8)
(m-6-9) edge (m-5-10) 
(m-6-11) edge (m-5-10);
\end{tikzpicture}
\end{equation*}
in which the backward pointing maps are all of the form
\[p_{k_{i,j}|k_{i+1,j}},\]
and the forward pointing maps are all of the form
\[j_{k_{i,j}|k_{i,j+1}}.\]
It's a trivial matter to see that this assignment defines a simplicial map
\[\widetilde{M}_S\colon\fromto{O(\Phi_S)}{A^{\eff}(\Fin)},\]
as desired.


\section{Localizing exact $\infty$-categories}

Now we apply this when $C=\categ{Exact}_{\infty}$, the $\infty$-category of exact $\infty$-categories. In particular, since exact $\infty$-categories form an $\infty$-category with direct sums, 
we may form, for any exact $\infty$-category $E$, the exact $\infty$-category $S^{-1}E$ via this filtered colimit. Since the multiplication by $k$ functor $k\colon\fromto{E}{E}$ induces the the multiplication by $k$ map $k\colon\fromto{K(E)}{K(E)}$, and since algebraic $K$-theory preserves filtered colimits, we deduce that
\[K(S^{-1}E)\simeq S^{-1}K(E).\]

The question, now, is whether our exact $\infty$-category $S^{-1}E$ is at all understandable. Happily, the answer is yes: we can identify $S^{-1}E$ with an $\infty$-category of certain graded objects, not quite of $E$, but of a natural enlargement thereof, where we might find suitably infinite objects for our analysis.
\begin{dfn} If $E$ is an essentially small exact $\infty$-category, then a \emph{large} object of $E$ is a functor $\fromto{E^{\op}}{\Top}$ that carries any zero object in $E$ to a terminal object and any admissible pushout/pullback square
\begin{equation*}
\begin{tikzpicture}[baseline]
\matrix(m)[matrix of math nodes,
row sep=4ex, column sep=4ex,
text height=1.5ex, text depth=0.25ex]
{X & Y \\
X' & Y' \\ };
\path[>=stealth,->,font=\scriptsize]
(m-1-1) edge[>->] node[above]{} (m-1-2)
edge[->>] node[left]{} (m-2-1)
(m-1-2) edge[->>] node[right]{} (m-2-2)
(m-2-1) edge[>->] node[below]{} (m-2-2);
\end{tikzpicture}
\end{equation*}
to a pullback square. We write $P_+(E)$ for the full subcategory of $\Fun(E^{\op},\Top)$ spanned by the large objects of $E$.
\end{dfn}

It is easy to check that $P_+(E)$ is a compactly generated, additive $\infty$-category, and that the Yoneda embedding of $E$ into $P_+(E)$ carries admissible pushout/pullback squares to squares that are both pushout and pullback squares. We may declare a morphism of $P_+(E)$ to be ingressive or egressive if and only if it is a filtered colimit of ingressive or egressive morphisms of $E$, respectively. With this structure, $P_+(E)$ is an exact $\infty$-category, and $j_+\colon\into{E}{P_+(E)}$ is exact.

Furthermore, $P_+(E)$ has the following universal property: for any additive, presentable $\infty$-category $D$, precomposition with $j_+$ defines an equivalence
\[
\equivto{\Fun^L(P_+(E),D)}{\Fun_{\Exact_{\infty}}(E,D).}
\]

\begin{exm} When $E$ is the ordinary category of finitely generated projective modules over a commutative ring $R$, then $P_+(E)$ is equivalent to the $\infty$-category $\categ{Ch}^+(R)$ of nonnegative chain complexes of $R$-modules.
\end{exm}

\begin{exm} More generally, when $E$ has its \emph{minimal} exact structure, so that the only ingressive morphisms are summand inclusions, the $\infty$-category $P_+(E)$ is the nonabelian derived $\infty$-category of $E$.
\end{exm}

\begin{exm} When $E$ is a stable $\infty$-category with its maximal exact structure, so that every morphism is ingressive, the $\infty$-category $P_+(E)$ is simply $\Ind(E)$.
\end{exm}

\begin{dfn} Suppose again $E$ an exact $\infty$-category and $S$ a set of primes. Then an \emph{$S$-divisible large object} of $E$ is an object of the (homotopy) limit of the functor
\[
P_+(E)[S]^{\vee}\colon\fromto{\Phi_S^{\op}}{\Cat_{\infty}}.
\]
We write $\categ{Div}_S(P_+(E))$ for this homotopy limit.

More concretely, an $S$-divisible object is a sequence of large objects
\[
\left\{X_i\right\}_{i\in\Phi_S}
\]
along with equivalences
\[
\rho_{i,j}\colon\equivto{X_i}{jX_{ij}}
\]
for any $i,j\in\Phi_S$, which fit together to give, for every $i_0,i_1,\dots,i_n\in\Phi_S$, an $n$-simplex
\[
X_{i_0}\equivto{}{} i_1X_{i_0i_1}\equivto{}{}\cdots\equivto{}{}i_1i_2\cdots i_nX_{i_0i_1\cdots i_n}
\]
of equivalences.
\end{dfn}

\begin{ntn} For any $m\in\Phi_S$, we have the projection
\[
\omega_m\colon\fromto{\categ{Div}_S(P_+(E))}{P_+(E)},
\]
given by evaluation at $m\in\Phi_S$, and we also have its left adjoint $\sigma_m$.

Given an object $V$ of $E$ and a natural number $m$, we may define an $S$-divisible large object
\[
\frac{V}{m}\coloneq\sigma_m(j_+(V)).
\]
We write $\categ{Div}_S(E)$ for the full subcategory of $\categ{Div}_S(P_+(E))$ spanned by the objects of the form $\frac{V}{m}$.
\end{ntn}

\begin{nul} Note that if $n=mk$ in $\Phi_S$, then 
\[
m\frac{V}{n}\simeq\frac{V}{k},
\]
justifying our notation.
\end{nul}


\begin{thm} Suppose $E$ an exact $\infty$-category and $S$ a set of primes. Then the exact $\infty$-category $S^{-1}E$ is equivalent to $\categ{Div}_SE$.
\begin{proof} The $\infty$-category $S^{-1}E$ is the colimit of the diagram $E[S]\colon\fromto{\Phi_S}{\Cat_{\infty}}$. We consider the embedding $\into{E}{P_+(E)}$, which is visibly functorial in $S$ and lands in the subcategory of compact objects. Hence the induced functor $\into{S^{-1}E}{S^{-1}P_+(E)}$ is fully faithful and exact, where $S^{-1}P_+(E)$ is computed in the $\infty$-category $\Pr^L$.

Now $S^{-1}P_+(E)$ is by definition the filtered colimit of $P_+(E)[S]$ computed in $\Pr^L$, which is in turn equivalent to the filtered limit of the adjoint diagram in $\Pr^R$, which is in turn the limit in $\Cat_{\infty}$. The adjoint diagram is clearly $P_+(E)[S]^{\vee}$, whence we find that $S^{-1}P_+(E)\simeq\categ{Div}_SP_+(E)$.

Now the essential image of the functor is spanned by those objects that lie in the image of an object $V$ of $E$ lying in some degree $m\in\Phi_S$. These are exactly the objects $\frac{V}{m}$ defined above.
\end{proof}
\end{thm}

\begin{exm} In the particular case in which $E$ is an idempotent complete stable $\infty$-category, the $\infty$-category $S^{-1}E\simeq\categ{Div}_S(E)$ is the full subcategory of $\categ{Div}_S(\Ind(E))$ spanned by the compact objects.
\end{exm}

\begin{rem} If $E$ is a symmetric monoidal exact $\infty$-category (i.e., an exact $\infty$-category whose underlying Waldhausen $\infty$-category is symmetric monoidal in the sense of \cite{MR3398729}), then one can show that $S^{-1}E$ is naturally an $E$-module, and the functors $\sigma_m\circ j_+\colon\fromto{E}{S^{-1}E}$ are $E$-module functors.
\end{rem}

\begin{wrn}\label{wrn:Slocal} We stress that $S^{-1}E$ will not in general be an $S$-local exact $\infty$-category. In fact, it is not hard to see that the only $S$-local exact $\infty$-category is $0$.
\end{wrn}




\section{Divisible objects as equivariant sheaves} In this section we will find a more geometric description of $S^{-1}C$ when $C$ is a presentable exact category, such as $P_+(E)$, and then we will cut the resulting large $\infty$-category back down to size. To begin, let us describe an action of the \emph{$S$-adic circle group} $\TT_S=S^{-1}\ZZ/\ZZ$ on the Cantor space $\Omega$.

\begin{ntn} For any prime number $p$, write
\[
\Omega_p\coloneq\Map(\NN,\langle p\rangle),
\]
equipped with the product topology. This is of course a Cantor space, as is the product
\[
\Omega_S\coloneq\prod_{p\in S}\Omega_p.
\]
(Of course $\Omega_p$ may be identified with the group $\ZZ_p$ of $p$-adic integers, but we won't use much of the abelian group structure.)

For any nonnegative integer $n$, we obtain a continuous map
\[
p^n\colon\fromto{\Omega_p}{\Omega_p},
\]
which carries $r$ to the map given by
\[
(p^nr)_i=\begin{cases} 0&\text{if }i\leq n;\\ r_{i-n}&\text{if }i>n.
\end{cases}
\]
(In other words, this is multiplication by $p^n$ in $\ZZ_p$.) For any product $m=\prod_{p\in S}p^{\nu_p(m)}$ of primes in $S$, we therefore obtain a continuous map
\[
m\colon\fromto{\Omega_S}{\Omega_S}.
\]
We write $m\Omega_S\subseteq\Omega_S$ for the image of this map, which is again a Cantor space. There is also a surjection $f_{p^n}\colon\fromto{\Omega_p}{p^n\Omega_p}$ given by
\[
f_{p^n}(r)_i=\begin{cases} 0&\text{if }i\leq n;\\ r_i&\text{if }i>n;
\end{cases}
\]
this extends to a surjection $f_{m}\colon\fromto{\Omega_S}{m\Omega_S}$ for any natural number $m$.
\end{ntn}

\begin{cnstr}\label{prp:Q/Z-action} Of course we have the free action of the cyclic group $C_p$ on $\langle p\rangle$, which clearly extends to a free action of $\TT_p$ on $\Omega_p$. Moreover, two elements $x,y\in\Omega_p$ lie in the same orbit if and only if $f_{p^n}(x)=f_{p^n}(y)$ for some nonnegative integer $n$.

These actions together provide an action of $\TT_S\cong\bigoplus_{p\in S}\TT_p$ on $\Omega_S$, and two elements $x,y\in\Omega_S$ lie in the same orbit if and only if $f_{m}(x)=f_{m}(y)$ for some natural number $m$.
\end{cnstr}

\begin{prp}\label{prp:rat-category}
    Let $C$ be an exact presentable $\infty$-category (e.g. $P_+(E)$ for an exact $\infty$-category $E$). Then there is an equivalence
    \[S^{-1}C \simeq \categ{Sh}^{\TT_S}_{C}(\Omega_S)\]
    where the right hand side is the $\infty$-category of $C$-valued $\TT_S$-equivariant sheaves on the space $\Omega_S$, with the $S$-adic circle group $\TT_S$ acting as above.
\end{prp}
\begin{proof}
The category $S^{-1}C$ is the colimit of a diagram $\Phi_S\to\Pr^L$. We can interpret the arrows appearing in this diagram as formed via a push-pull construction
\[C\xrightarrow{\pi^*} \categ{Sh}_C\left(\langle n\rangle\right) \xrightarrow{\pi_*} C\]
where $\pi\colon\fromto{\langle n\rangle}{\langle 1\rangle}$ is the projection. But we can decouple the pullback and the pushforward by employing \S \ref{app:constructESandES} to define a factorization of $\Phi_S\to \Pr^L$ through a functor $O(\Phi_S)=\Fun(\Delta^1,\Phi_S)\to \Pr^L$ that carries each object $(m|n)$ of $O(\Phi_S)$ to the $\infty$-category
\[
\categ{Sh}_C\left(\left\langle\frac{n}{m}\right\rangle\right).
\]
Precisely, we compose $\widetilde{M}_S$ with the unique functor $\categ{Sh}_C\colon\fromto{A^{\eff}(\Fin)}{\Pr^L}$ that preserves finite products and carries $\langle 1\rangle$ to $C$ (with the direct sum symmetric monoidal structure).

Since $\Phi_S$ is a filtered category, the inclusion $\Phi_S\to O(\Phi_S)$ is cofinal and we can compute
\begin{equation*}
    S^{-1}C \coloneq \colim_{m\in \Phi_S} C\simeq\colim_{(m|n)\in O(\Phi_S)} \categ{Sh}_C\left(\left\langle\frac{n}{m}\right\rangle\right)\simeq \colim_{m\in \Phi_S} \colim_{n\in\Phi_S,\ m|n} \categ{Sh}_C\left(\left\langle\frac{n}{m}\right\rangle\right)
\end{equation*}
where in the last equality we have used that the projection $O(\Phi_S)\to \Phi_S$ sending $(m|n)$ to $m$ is a cocartesian fibration and so we can compute colimits fiberwise. But, since colimits in $\Pr^L$ can be computed as limits in $\Pr^R$, we have for any fixed $m\in\Phi_S$,
\[\colim_{n\in\Phi_S,\ m|n} \categ{Sh}_C\left(\left\langle\frac{n}{m}\right\rangle\right) = \lim_{n\in\Phi_S,\ m|n} \categ{Sh}_C\left(\left\langle\frac{n}{m}\right\rangle\right)=\categ{Sh}_C\left(m\Omega_S\right).\]
Here, the final identification follows from the fact that the $\infty$-category of sheaves on the lattice of clopen sets $m\Omega_S$ (i.e., the union of the lattices of subsets of $\langle\frac{n}{m}\rangle$ as $n$ varies through $\Phi_S$) is equivalent to the $\infty$-category of sheaves on the topological space $m\Omega_S$, because clopen sets form a basis that is closed under finite intersections.

So we have shown that
\[S^{-1}C\simeq \colim_{m\in \Phi_S} \categ{Sh}_C\left(m\Omega_S\right),\]
where the maps in the diagram are given by the pushforward along the projection \[j_{m|n}\colon m\Omega_S\to n\Omega_S.\] But colimits in $\Pr^L$ can be computed as limits in $\Pr^R$ after replacing all the functors with their right adjoints. Since $j_{m|n}$ is étale and proper, the right adjoint of the pushforward is the pullback. hence we can write
\[S^{-1}C \simeq \lim_{m\in \Phi_S^{\op}} \categ{Sh}_C\left(m\Omega_S\right)\,.\]

Now we observe that the map $j_{1|m}\colon\Omega_S\to m\Omega_S$ is the surjection $f_m$ above. In particular we can write
\[S^{-1}C\simeq\lim_{m\in\Phi_S^{\op}}\categ{Sh}_C(\Omega_S/R_m)\simeq \lim_{m\in\Phi_S^{\op}}\lim_{\Delta^\op} \categ{Sh}_C\left(R_m\times_{\Omega_S} \cdots \times_{\Omega_S}R_m\right)\,,\]
where $R_m$ is the equivalence relation given by
\[R_m=\left\{(x,y)\in\Omega_S\times\Omega_S\ \middle|\ f_m(x)=f_m(y)\right\}\,,\]
and we conclude that
\[
S^{-1}C\simeq\lim_{\Delta^\op} \categ{Sh}_C\left(R\times_{\Omega_S} \cdots \times_{\Omega_S}R\right)\,,
\]
where $R=\colim_{m\in\Phi_S}R_m$.
Finally, by Cnstr. \ref{prp:Q/Z-action}, the equivalence relation $R$ is exactly the equivalence relation induced on $\Omega_S$ by the action of $\TT_S$. So
\[S^{-1}C \simeq\categ{Sh}_C\left(\Omega_S\right)^{h\TT_S}\,,\]
as desired.
\end{proof}

\begin{rem} A simple analysis of this proof shows that if $C$ is a presentably symmetric monoidal exact $\infty$-category, then the equivalence $S^{-1}C\simeq \categ{Sh}_C^{\TT_S}(\Omega_S)$ is an equivalence of $C$-modules.
\end{rem}


\begin{nul} Note that since $\Omega_S$ is a compact Hausdorff space of finite covering dimension, it follows that the corresponding $\infty$-topos is hypercomplete. This ensures that equivalences in $\categ{Sh}_C(\Omega_S)$ and $\categ{Sh}_C^{\TT_S}(\Omega_S)$ can be detected on stalks.
\end{nul}

\begin{nul} Of course we wish to apply this to the case in which $C=P_+(E)$ for some exact $\infty$-category $E$. The full subcategory $S^{-1}E\subset S^{-1}P_+(E)$ can be identified with a full subcategory
\[\categ{Sh}_{P_+(E)}^{\TT_S}(\Omega_S)^{\textit{small}}\subseteq\categ{Sh}_{P_+(E)}^{\TT_S}(\Omega_S).\]

The objects $\frac{V}{m}$ of $\categ{Sh}_{P_+(E)}^{\TT_S}(\Omega_S)^{\textit{small}}$ can be described as follows. Form the constant sheaf $V$ on $\langle m\rangle$ with the obvious $C_m$ action; call the result $V$ again. Now $\frac{V}{m}$ is the induced $\TT_S$-equivariant sheaf
\[
\TT_S\times_{C_m}V\cong\bigoplus_{g\in\TT_S/C_m}g^{\star}V
\]
on $\Omega_S$.

Now if $E$ is an idempotent-complete stable $\infty$-category, then $\categ{Sh}_{P_+(E)}^{\TT_S}(\Omega_S)^{\textit{small}}$ is the full subcategory of $\categ{Sh}_{\Ind E}^{\TT_S}(\Omega_S)$ spanned by the compact objects.

If $E$ is a symmetric monoidal exact $\infty$-category, then one can show that the $P_+(E)$-module equivalence $S^{-1}P_+(E)\simeq\categ{Sh}_{P_+(E)}^{\TT_S}(\Omega_S)$ restricts to an $E$-module equivalence $S^{-1}E\simeq \categ{Sh}_{P_+(E)}^{\TT_S}(\Omega_S)^{\textit{small}}$.
\end{nul}

We now turn our attention to the $G$-theory of a quasi-compact quasi-separated scheme $X$. (Everything will also work in the derived or spectral settings with small modifications that are best left to the reader.) Following Illusie, one defines the $\infty$-category $\categ{Coh}(X)\subset\categ{QCoh}(X)$ of \emph{coherent complexes} on $X$ as follows:
\begin{enumerate}[(1)]
\item If $X=\Spec A$ is an affine scheme, then $\categ{Coh}(X)$ is defined as the full subcategory of the derived $\infty$-category $\DD(A)$ spanned by those bounded complexes of $A$-modules $M$ such that for any filtered diagram $\{N_{\alpha}\}_{\alpha\in\Lambda}$ of $A$-modules, and any integer $n$, the natural map
\[
\fromto{\colim_{\alpha\in\Lambda}\Map(M,N_{\alpha}[n])}{\Map(M,\colim_{\alpha\in\Lambda}N_{\alpha}[n])}
\]
is an equivalence.
\item In general, an object of $\categ{QCoh}(X)$ belongs to the subcategory $\categ{Coh}(X)$ if and only if its restriction to every affine open subscheme $U\subset X$ belongs to $\categ{Coh}(U)$. We set
\[
\categ{IndCoh}(X)\coloneq\Ind\categ{Coh}(X).
\]
\end{enumerate}
Recall that the $G$-theory of $X$ is defined by
\[
G(X)\coloneq K(\categ{Coh}(X)).
\]

Now recall that $\Omega_S$ can be seen as an affine scheme (precisely as the spectrum of the ring of locally constant $\ZZ$-valued functions on $\Omega_S$). 
Since
\[
\categ{Sh}_{\categ{IndCoh}(X)}(\Omega_S)\simeq\categ{IndCoh}(X\times\Omega_S)
\]
we can express Pr. \ref{prp:rat-category} in a different way:
\begin{prp} Let $X$ be a quasi-compact quasi-separated scheme. There is an equivalence of stable presentable $\infty$-categories
\[
S^{-1}\categ{IndCoh}(X)\simeq\categ{IndCoh}(X\times\Omega_S)^{h\TT_S}\,.
\]
\end{prp}
\noindent Following Gaitsgory \cite{Gaitsgory}, we may extend the definition of $\categ{IndCoh}$ to more general objects by stipulating that the functor
$\goesto{X}{\categ{IndCoh}(X)}$, $\goesto{f}{f^!}$, transform colimits into limits.
 The quotient algebraic space
\[\left[\left(X\times \Omega_S\right)/\TT_S\right]\simeq X\times \left[\Omega_S/\TT_S\right]\]
can be expressed as a colimit of schemes
\[
\colim_{m\in\Phi_S} \left(X\times \Omega_S\right)/C_m
\]
in which all maps are finite étale. Since $f^!=f^*$ for such maps $f$, we obtain
\[
S^{-1}\categ{IndCoh}(X) \simeq \categ{IndCoh}\left(\left[\left(X\times \Omega_S\right)/\TT_S\right]\right).
\]


As $S^{-1}\categ{IndCoh}(X)$ is furthermore compactly generated \cite[Pr. 5.5.7.6]{HTT}, it is sensible to define $\categ{Coh}\left(\left[\left(X\times \Omega_S\right)/\TT_S\right]\right)$ as the full stable subcategory of the $\infty$-category $\categ{IndCoh}\left(\left[\left(X\times \Omega_S\right)/\TT_S\right]\right)$ spanned by the compact objects. Consequently, the proposition above induces an identification
\[
S^{-1}\categ{Coh}(X)\simeq\categ{Coh}^{\TT_S}(X\times\Omega_S)=\categ{Coh}\left(\left[\left(X\times \Omega_S\right)/\TT_S\right]\right).
\]
We thus obtain the desired identification of spectra (and even $K(X)$-modules)
\[
S^{-1}G(X)\simeq G^{\TT_S}(X\times\Omega_S).
\]
In particular, when $X=\Spec A$, then one has
\[
S^{-1}G(A)\simeq G^{\TT_S}(C(\Omega_S,A)),
\]
where $C$ denotes the ring of locally constant functions.

\begin{rem} We caution that the algebraic space $\left[\left(X\times \Omega_S\right)/\TT_S\right]$ is not perfect: compact objects such as $\frac{O_X}{1}$ are not dualizable in the symmetric monoidal $\infty$-category $\categ{QCoh}\left(\left[\left(X\times \Omega_S\right)/\TT_S\right]\right)$, and conversely the unit object is not compact. Hence, we cannot simply replace $G$-theory by $K$-theory in the above formulas.
\end{rem}

\bibliographystyle{plain}
\bibliography{ratK}

\end{document}